\documentclass[nolinenbrs]{dmgt-silent}
\newif\ifdraft
\draftfalse

\usepackage[T1]{fontenc}
\usepackage[utf8]{inputenc}

\usepackage[centertags,reqno]{amsmath}
\usepackage{cite}
\usepackage{fancyhdr}
\usepackage{mathtools}

\newtheorem{definition}[theorem]{Definition}
\newtheorem{proposition}[theorem]{Proposition}
\newtheorem{corollary}[theorem]{Corollary}

\DeclareMathOperator*{\adj}{adj}

\newcommand{\dotcup}{\mathrel{\dot{\cup}}}

\renewcommand{\vec}[1]{\ensuremath{\mathbf{#1}}}

\newcommand{\field}[1]{\ensuremath{\mathbb{#1}}}

\newcommand{\comp}[1]{\ensuremath{\overline{#1}}}
\newcommand{\cp}[1]{\ensuremath{\mathop{\phi}\left(#1,\lambda\right)}}

\newcommand{\inserttab}[4]{%
   \begin{table}[!p]
      \begin{minipage}{\columnwidth}
         \renewcommand\footnoterule{}
         \caption{#3}
         \label{#1}
         \centering
         \renewcommand{\arraystretch}{1.15}
         \begin{tabular}{#2}
            \hline
            #4
            \hline
         \end{tabular}
      \end{minipage}
   \end{table}}

\newcommand\gitRevision{6d27a93}
\newcommand\gitAuthorDate{Wed Jan 29 13:45:34 2020 +0100}

\def\tfn{\footnote{Manuscript submitted March 28, 2019;
revised August 2, 2019;
accepted August 22, 2019 for publication in Discussiones Mathematicae Graph Theory.\\
Dedicated to the memory of Slobodan K. Simi\'c.}}

\title{On the Displacement of Eigenvalues when Removing a Twin Vertex\tfn}
[On the Displacement of Eigenvalues]

\keywords{%
   eigenvalues,
   perturbations,
   duplicate and co-duplicate vertices,
   threshold graph,
   nested split graph}

\classnbr{05C50} 

\newauthor{Johann~A.~Briffa}{J.A.Briffa}
   {Department of Communications and Computer Engineering\\
   University of Malta, Msida MSD~2080, Malta}
   [johann.briffa@um.edu.mt]
\newauthor{Irene~Sciriha}{I.Sciriha}
   {Department of Mathematics\\
   University of Malta, Msida MSD~2080, Malta}
   [irene.sciriha-aquilina@um.edu.mt]

\begin{document}

\thispagestyle{fancy}
\fancyhf{}
\fancyfoot[L]{\footnotesize\itshape r.\gitRevision}
\fancyfoot[C]{}
\fancyfoot[R]{\footnotesize\itshape \gitAuthorDate}
\renewcommand{\headrulewidth}{0pt}
\renewcommand{\footrulewidth}{0pt}

\begin{abstract}
Twin vertices of a graph have the same open neighbourhood.
If they are not adjacent, then they are called duplicates and contribute the
eigenvalue zero to the adjacency matrix.
Otherwise they are termed co-duplicates, when they contribute $-1$ as an
eigenvalue  of the adjacency matrix.
On removing a twin vertex from a graph, the spectrum of the adjacency matrix
does not only lose the eigenvalue $0$ or $-1$.
The perturbation sends a rippling effect to the spectrum.
The simple eigenvalues are displaced.
We obtain a closed formula for the characteristic polynomial of a graph with
twin vertices in terms of two polynomials associated with the perturbed graph.
These are used to obtain estimates of the displacements in the spectrum
caused by the perturbation.
\end{abstract}

\clearpage

\section{Introduction}

We limit ourselves to simple connected graphs, that is graphs with no multiple
edges or loops.
A graph $G(V,E)$ has a vertex set $V=\{1,2,\ldots,n\}$ and an edge set $E$
whose elements are distinct pairs of vertices of $V$.
The set $\comp{E}$ of non-edges of $G$ are those pairs of distinct vertices
not in $ E$.
The complement $\comp{G}(V,\comp{E})$ of $G$ has the same vertex set as $G$
and edge set $\comp{E}$.
Twin vertices are either duplicate or co-duplicate.
Two vertices are called duplicate if they are non-adjacent and have the
same neighbours.
A pair of co-duplicate vertices in a graph $G$ are adjacent, and they are
duplicate vertices in the complement $\comp{G}$.

Let $\vec{A}(G)$, also written as $\vec{A} = ( a_{i,j} )$, be the adjacency
matrix of $G$ with $a_{i,j}=1$ if the vertices $i,j$ are adjacent and zero
otherwise.
The eigenvalues of $\vec{A}$ are referred to as the eigenvalues of $G$
and form the spectrum of $G$.
If $G$ has a pair of duplicate vertices, then the corresponding rows (and
columns) in $\vec{A}$ are the same.
This means that $\vec{A}$ has the eigenvalue zero.
In the case where $G$ has two co-duplicate vertices, the corresponding
rows and columns are the same except for the two entries defining the
edge between them.
This means that $-1$ is in the spectrum of $G$.
In both cases the associated eigenvector has two non-zero entries.


Unlike what one may assume, removing a twin vertex does not just remove
the eigenvalue $0$ or $-1$ in the respective cases, while preserving
the rest of the spectrum.
Indeed, we investigate the shift in eigenvalues on removing a twin vertex.
To calculate the new eigenvalues after removing a twin vertex, one has to
perform the computation on the adjacency matrix of the new graph, ignoring
any information known about the original graph.
In this work, we provide ways to directly calculate estimates for the changes
in eigenvalues, as a difference from those of the original graph.
We also give an explicit expression for the change in the characteristic
polynomial due to the removal of a twin vertex.


While to our knowledge this specific problem has not been treated before,
the literature on spectral graph theory contains a number of related works.
In the 1950s, Heilbronner derived the characteristic polynomial of a perturbed
graph from that of the parent graph.
He determined explicitly the eigenvalues of the subgraph on deleting a
vertex from a graph, contributing to the study of molecular orbitals
\cite{heilbronner1953kompositions, heilbronner1954molecular,
heilbronner1954graphisches, heilbronner1962graphentheoretischen,
heilbronner1979some}.
Later, in the literature, one finds expressions for the characteristic
polynomial of an arbitrary graph, of graphs with particular geometric
properties and of perturbed graphs also in the work of Schwenk
\cite{schwenk1974} and Rosenfeld \cite{rosenfeld2013}.

The well known Cauchy inequalities, involving the eigenvalues of a real
symmetric matrix and a principal submatrix, are referred to as the interlacing
theorem in spectral graph theory \cite{schwenk1974}.
The theorem states that exactly one root of the characteristic polynomial of
a vertex deleted graph lies between two successive eigenvalues of the parent
graph.
It was the subject of many studies by the pioneers in the theory of the
matrices that encode the structure of a graph.
Its application unlocked many remarkable latent properties of classes
of graphs.
Interlacing is the main tool used by Thüne in his PhD thesis \cite{thune1979}
to determine certain substructures in graphs. 
Later, Haemers produced a survey \cite{haemers1995} of the various kinds
of applications of eigenvalue interlacing to complement his doctoral thesis
\cite{haemers1980}.
More recently, Lovász emphasised the importance of interlacing and gave
a summary of the main results on the eigenvalues of matrices of a graph
\cite{lovasz2007}.
In a collaboration with Simić, Marino \emph{et al.}\ used the interlacing
theorem to explore the properties of line graphs of trees with twin vertices
deleted \cite{marino2006}.
Sciriha \emph{et al.}\ studied classes of graphs that showed the largest
possible change in the multiplicity of the eigenvalue zero as a consequence
of interlacing \cite{sdbfp2013}.
An independent set is a subset of the vertices such that no two of the
vertices are adjacent.
One of the graph invariants that is widely studied in combinatorics is the
independence number, that is the size of the largest independent set.
Rowlinson, a main exponent of graph eigenvalues, finds bounds on the
independence number of a graph basing his arguments on interlacing
\cite{rowlinson2007}.

Sciriha and Farrugia consider the threshold graph which is a split graph in
which the vertex set is partitioned into an independent set and a clique,
which is a subset of the vertices that induces a complete subgraph in which
every two vertices are adjacent \cite{sciriha2012spectrum}.
The independent set may contain duplicates and the clique may
contain co-duplicates.
Mohammadian and Trevisan show that there are no eigenvalues of the adjacency
matrix of a threshold graph between $0$ and $-1$, which in threshold graphs
are contributed only by twin vertices \cite{mohammadian2016}.

The interlacing theorem applies to all the real symmetric matrices encoding
$G$, including the Laplacian matrix.
When considering the Laplacian, So proved that only one of its eigenvalues
is displaced when an edge is added between two duplicate vertices to produce
a co-duplicate \cite{so1999}.
The rest of the eigenvalues remain unchanged.


The interlacing theorem provides rough bounds for the displacement of the
eigenvalues of $G$ when a vertex is deleted.
Our objective is to obtain better estimates within these bounds.
To this end, relations of $\cp{\vec{A}(G)}$ to polynomials of other
submatrices of $\vec{A}$ are obtained.
These results would be of interest in any application where the displacement in
eigenvalues is of greater interest than the eigenvalues of the modified graph.


The rest of the paper is organised as follows.
In Section~\ref{sec:formulae}, we apply similarity operations on the
adjacency matrix of $G$, so that eigenvalues are preserved, to yield a
matrix whose characteristic polynomial is easily expressed in terms of
those of subgraphs of $G$.
In Section~\ref{sec:displacement}, we show how the expressions obtained
enable the computation of estimates of the displacement of the eigenvalues
of the adjacency matrix on removing a twin vertex.
Finally, we give examples of computing the estimates of the displacement
of the spectrum on removing a twin vertex from a nested split graph
in Appendix~\ref{sec:examples_nsg}, and from a general graph in
Appendix~\ref{sec:examples_nut}.


\section{Effect on the Characteristic Polynomial on Removing a Twin Vertex}
\label{sec:formulae}

To obtain the eigenvalues of a matrix $\vec{M}$, it suffices to determine
the roots of its characteristic polynomial $\cp{\vec{M}}$.
If $\vec{M}$ is known to be real and symmetric, then its algebraic
properties allow alternative methods of computation with possibly lower
complexity.
The Jacobi-Givens method \cite{froberg} employs rotation of two axes
of $\mathbb{R}^n$ to introduce zero entries in a row of $\vec{M}$ via a
similarity operation and therefore without altering the eigenvalues.
The new form of the matrix allows the characteristic polynomials of $\vec{M}$
and of other principal submatrices of $\vec{M}$ to be easily related.

\begin{definition}[(Adjacency matrix)]
The adjacency matrix $\vec{A}$ of a graph $G$ of order $n$, where the two first
labelled vertices $v_1$,$v_2$ are twin vertices, can be written as
\begin{equation}
\label{eqn:A}
\vec{A} = \left(
   \begin{array}{cc|c}
   0 & a & \vec{b}^{\top} \\
   a & 0 & \vec{b}^{\top} \\
   \hline
   \vec{b} & \vec{b} & \vec{C}
   \end{array}
\right)
\text{,}
\end{equation}
where $\vec{C}=\vec{A}(G_{-v_1-v_2})$ is the adjacency matrix of the subgraph
$G_{-v_1-v_2}$ of $G$, obtained from $G$ by removing vertices $v_1$,$v_2$ and
the edges incident to them.
The entry $a$ is $0$ for duplicate and $1$ for co-duplicate vertices.
\end{definition}

\begin{proposition}
\label{prop:givens}
The adjacency matrix $\vec{A}$ is similar to the simpler matrix
\begin{equation}
\vec{A'}
   = \left(
      \begin{array}{cc|c}
      a & 0  & \sqrt{2}\vec{b}^{\top} \\
      0 & -a & \vec{0}^{\top} \\
      \hline
      \sqrt{2}\vec{b} & \vec{0} & \vec{C}
      \end{array}
   \right)
\text{.}
\end{equation}
\end{proposition}

\begin{proof}
We use the Jacobi-Givens method to find a matrix $\vec{P}$ such that
$\vec{A'} = \vec{P^{-1}}\vec{A}\vec{P}$.
Since twin vertices have the same open neighbourhood, a rotation by
$\frac{\pi}{4}$ of the corresponding axes in $\field{R}^n$ is required.
This is achieved by using
\begin{equation}
\vec{P} = \left(
   \begin{array}{cc|c}
   \frac{1}{\sqrt{2}} & -\frac{1}{\sqrt{2}} & \vec{0} \\
   \frac{1}{\sqrt{2}} &  \frac{1}{\sqrt{2}} & \vec{0} \\
   \hline
   \vec{0} & \vec{0} & \vec{I}
   \end{array}
\right)
\text{,}
\end{equation}
where $\vec{I}$ is the identity matrix.
\end{proof}

\begin{corollary}
\label{cor:givens}
\begin{equation}
\cp{\vec{A'}}
   = \det(\lambda \vec{I} - \vec{P^{-1}}\vec{A}\vec{P})
   = \left|
      \begin{array}{cc|c}
      \lambda - a & 0 & -\sqrt{2}\vec{b}^{\top} \\
      0 & \lambda + a & \vec{0}^{\top} \\
      \hline
      -\sqrt{2}\vec{b} & \vec{0} & \lambda \vec{I} - \vec{C}
      \end{array}
   \right|
\text{.}
\end{equation}
\end{corollary}

\begin{proposition}
The characteristic polynomial of a graph $\vec{G}$ with adjacency matrix
$\vec{A}$ having a pair of twin vertices is
\begin{equation}
\label{eqn:PHIA}
\cp{\vec{A}}
   = (\lambda^2 -a^2) \cp{\vec{C}}
      - 2 (\lambda + a) \vec{b}^\top \adj(\lambda \vec{I} - \vec{C}) \vec{b}
\text{,}
\end{equation}
where the adjugate $\adj(\lambda \vec{I} - \vec{C})$ is equivalent to the
expression
\begin{equation*}
(\lambda \vec{I} - \vec{C})^{-1} \det(\lambda \vec{I} - \vec{C})
\text{,}
\end{equation*}
for non-singular $\lambda \vec{I} - \vec{C}$.
\end{proposition}

\begin{proof}
Using $\cp{\vec{A'}}$ from Corollary~\ref{cor:givens} we can express the
characteristic polynomial of $\vec{A}$ as
\begin{equation}
\cp{\vec{A}}
   = \det(\lambda \vec{I} -\vec{A})
   = \det(\lambda \vec{I} -\vec{A'})
   = (\lambda + a) \left|
      \begin{array}{c|c}
      \lambda - a & -\sqrt{2}\vec{b}^{\top} \\
      \hline
      -\sqrt{2}\vec{b} & \lambda \vec{I} - \vec{C}
      \end{array}
   \right|
\text{,}
\end{equation}
written as $(\lambda + a) \det(\vec{M})$.
Expanding this expression in terms of the Schur complement
$\vec{M} | (\lambda \vec{I} - \vec{C})$ of $\vec{M}$,
\begin{align}
\cp{\vec{A}}
   &= (\lambda + a) \cp{\vec{C}} \det(\vec{M}|(\lambda \vec{I} - \vec{C})) \\
\cp{\vec{A}}
   &= (\lambda + a) \cp{\vec{C}} \left[
         (\lambda - a) - 2 \vec{b}^\top (\lambda \vec{I} - \vec{C})^{-1} \vec{b}
      \right]
\text{.}
\end{align}
The result follows immediately.
\end{proof}


\begin{lemma}
If $v_1$ is a twin vertex of $G$, then the characteristic polynomial
of the subgraph $G_{-v_1}$, obtained from $G$ by deleting vertex $v_1$,
is given by
\begin{equation}
\label{eqn:PHIAv}
\cp{\vec{A}(G_{-v_1})}
   = \lambda \cp{\vec{C}}
   - \vec{b}^\top \adj(\lambda \vec{I} - \vec{C}) \vec{b}
\text{.}
\end{equation}
\end{lemma}

\begin{proof}
Observe that
\begin{equation}
\cp{\vec{A}(G_{-v_1})}
   = \left|
      \begin{array}{c|c}
      \lambda & -\vec{b}^{\top} \\
      \hline
      -\vec{b} & \lambda \vec{I} - \vec{C}
      \end{array}
   \right|
\text{.}
\end{equation}
The result follows using the Schur complement expansion.
\end{proof}


Next, we obtain relations of $\cp{\vec{A}(G)}$ to polynomials of other submatrices
of $\vec{A}$.

\begin{definition}
\label{def:adjugate}
Let
$\adj(\lambda \vec{I} - \vec{A})
   = \left(
   h_{\ell,k}
   \right)_{n \times n}$
so that $h_{\ell,k}$ denotes the entry in row $\ell$ and column $k$ of the adjugate
$\adj(\lambda \vec{I} - \vec{A})$.
\end{definition}

\begin{lemma}
\label{lem:badjb}
Let $v_1$ and $v_2$ be twin vertices, and $\vec{C} = \vec{A}(G_{-v_1-v_2})$,
then
\begin{equation}
\label{eqn:h12}
h_{1,2}
   = a \cp{\vec{C}} + \vec{b}^\top \adj(\lambda \vec{I} - \vec{C}) \vec{b}
\text{.}
\end{equation}
\end{lemma}

\begin{proof}
The matrix $\adj(\lambda \vec{I} - \vec{A})$ is real and symmetric for
real $\lambda$.
So
\begin{equation}
h_{1,2}
   = h_{2,1}
   = -\left|
      \begin{array}{c|c}
      -a & -\vec{b}^\top \\
      \hline
      -\vec{b} & \lambda \vec{I} - \vec{C}
      \end{array}
   \right|
\text{.}
\end{equation}
The Schur complement expansion of the determinant, gives
\begin{align}
h_{1,2}
   & = -\cp{\vec{C}} \left[
      -a - \vec{b}^\top (\lambda \vec{I} - \vec{C})^{-1} \vec{b}
      \right] \\
   & = a \cp{\vec{C}} + \vec{b}^\top \adj(\lambda \vec{I} - \vec{C}) \vec{b}
\text{.}
\end{align}
\end{proof}


The characteristic polynomial of $\vec{A}$ in \eqref{eqn:A} can also be
expressed in terms of two determinants.

\begin{theorem}
\label{thm:G12}
Let the first two labelled vertices $v_1$  and $v_2$ be twin vertices.
Then
\begin{equation}
\cp{\vec{A}(G)}
   = (\lambda + a) \left[
         \cp{\vec{A}(G_{-v_1})}
         - h_{1,2}
      \right]
\text{.}
\end{equation}
\end{theorem}

\begin{proof}
Eliminating $\vec{b}^\top \adj(\lambda \vec{I} - \vec{C}) \vec{b}$
from~\eqref{eqn:PHIA} and \eqref{eqn:h12} we obtain
\begin{align}
\cp{\vec{A}(G)}
   & = (\lambda^2 - a^2) \cp{\vec{C}}
      - 2 (\lambda + a) \left[
         h_{1,2}
         - a \cp{\vec{C}}
      \right] \\
   \label{eqn:PhiA_simp}
   & = (\lambda + a)^2 \cp{\vec{C}}
      - 2 (\lambda + a) h_{1,2}
\end{align}
Similarly, eliminating $\vec{b}^\top \adj(\lambda \vec{I} - \vec{C}) \vec{b}$
from~\eqref{eqn:PHIAv} and \eqref{eqn:h12} we obtain
\begin{align}
\cp{\vec{A}(G_{-v_1})}
   & = \lambda \cp{\vec{C}} - \left[
         h_{1,2}
         - a \cp{\vec{C}}
      \right] \\
   \label{eqn:PhiAv_simp}
   & = (\lambda + a) \cp{\vec{C}} - h_{1,2}
\end{align}
Finally, eliminating $\cp{\vec{C}}$ from~\eqref{eqn:PhiA_simp}
and~\eqref{eqn:PhiAv_simp} completes the proof.
\end{proof}


\begin{lemma}
\label{lem:permutation}
Pre-multiplying a matrix
$\vec{M} = \left( m_{i,j} \right)_{n \times n}$
by the permutation matrix
\begin{equation*}
\vec{E}_{\ell,1} = \left(
   \begin{array}{c|c|c}
   0       & 1 & 0 \\
   \hline
   \vec{I}_{(\ell-1) \times (\ell-1)} & 0 & 0 \\
   \hline
   0       & 0 & \vec{I}_{(n-\ell) \times (n-\ell)} \\
   \end{array}
   \right)
\end{equation*}
gives $\vec{M'} = \left( m'_{i,j} \right)_{n \times n}$
with row $\ell$ of $\vec{M}$ in row 1 of $\vec{M'}$; that is the
entries of $\vec{M'}$ are given by
\begin{equation}
m'_{j,k} =
   \begin{cases}
      m_{\ell,k} & j = 1 \\
      m_{j-1,k} & 1 < j \leq \ell \\
      m_{j,k} & \text{otherwise}
   \end{cases}
\end{equation}
\end{lemma}

The effect of pre-multiplying $\vec{M}$ by $\vec{E}_{\ell,1}$ is to move
row $\ell$ of $\vec{M}$ to row 1 of $\vec{M'}$, shifting rows 1 to $\ell-1$
of $\vec{M}$ by one.
Post-multiplying $\vec{M}$ by the transpose of $\vec{E}_{\ell,1}$ has the
same effect on the columns.

\begin{proposition}
\label{prop:shift}
The matrix \vec{M''} is obtained by moving row $\ell$ and column $\ell$ of
$\vec{M}$ to the first row and first column, using
\begin{equation}
\vec{M''} = \vec{E}_{\ell,1} \vec{M} \vec{E}_{\ell,1}^\top
\end{equation}
\end{proposition}

The determinant of the product of two square matrices is the product of the
separate determinants.
Since $\vec{E}_{\ell,1}^\top = \vec{E}_{\ell,1}^{-1}$, the next result
follows immediately.

\begin{corollary}
\begin{equation}
\det(\vec{M}) = \det(\vec{M''})
\end{equation}
\end{corollary}

Recall that entry $\ell,k$ of the adjugate of a matrix is $h_{\ell,k}$,
the $\ell,k$ co-factor of the matrix.

\begin{proposition}
\label{prop:hlk}
\begin{equation}
h_{\ell,k} = (-1)^{\ell+k}
   \left|
   \begin{array}{c|c}
      -a_{\ell,k} & -\vec{b}_\ell^\top \\
      \hline
      -\vec{b}_k & \lambda \vec{I} - \vec{B} \\
   \end{array}
   \right|
\end{equation}
where $\vec{B}$ is obtained from $\vec{A}$ by deleting rows and
columns $\ell$ and $k$, $\ell \neq k$.
\end{proposition}

\begin{proof}
This follows immediately from Definition~\ref{def:adjugate}.
\end{proof}


Applying Proposition~\ref{prop:shift}, Theorem~\ref{thm:G12} can be
generalized to:

\begin{theorem}
\label{thm:Glk}
Let  $v_\ell$  and $v_k$ be twin vertices.
Then
\begin{equation}
\cp{\vec{A}(G)}
   = (\lambda + a_{\ell,k}) \left[
         \cp{\vec{A}(G_{-v_\ell})}
         - h_{\ell,k}
      \right]
\text{.}
\end{equation}
\end{theorem}


An alternative perspective is that we can obtain the characteristic polynomial
of the graph with a twin vertex removed.

\begin{corollary}
Let  $v_\ell$  and $v_k$ be twin vertices.
Then
\begin{equation}
\label{eqn:cp_change}
\cp{\vec{A}(G_{-v_\ell})}
   = \frac{ \cp{\vec{A}(G)} }{ \lambda + a_{\ell,k} }
      + h_{\ell,k}
\text{.}
\end{equation}
\end{corollary}


\section{Estimating the Displacement of Eigenvalues}
\label{sec:displacement}

In this section, the relation \eqref{eqn:cp_change} is used to obtain first
order and second order estimates of the displacement of eigenvalues on
deleting a twin vertex.
Define
\begin{equation}
f(\lambda)
   = \frac{ \cp{\vec{A}(G)} }{ \lambda + a_{\ell,k} }
      + h_{\ell,k}(\lambda)
\text{,}
\end{equation}
such that $\cp{\vec{A}(G_{-v_\ell})} = f(\lambda)$, which is a polynomial
in $\lambda$.
Now, we can express $f(\lambda)$ using the Taylor series
\begin{equation}
f(\lambda) =
   f(\lambda_0)
   + \frac{f'(\lambda_0)}{1!} (\lambda_0 - \lambda)
   + \frac{f''(\lambda_0)}{2!} (\lambda_0 - \lambda)^2
   + \cdots
\text{.}
\end{equation}
Choosing $\lambda_0$ to be a root of \cp{\vec{A}(G)} gives us an expression
in terms of $\delta = \lambda_0 - \lambda$, or the displacement from the
eigenvalue $\lambda_0$ when $f(\lambda) = 0$.
For a first order approximation, we truncate the Taylor series to the first
power of $\delta$, obtaining
\begin{align}
0 &=
   f(\lambda_0)
   + \delta f'(\lambda_0) \\
\delta &= -\frac{f(\lambda_0)}{f'(\lambda_0)}
\text{.}
\end{align}
Similarly, a second order approximation can be obtained by solving the
quadratic equation
\begin{equation}
\label{eqn:quadratic}
0 =
   f(\lambda_0)
   + \delta f'(\lambda_0)
   + \delta^2 \frac{f''(\lambda_0)}{2}
\text{.}
\end{equation}
The displacement depends on the mapping of the eigenvalues of the original
graph to those of the vertex-deleted subgraph.
This mapping is uniquely determined by retaining the order of eigenvalues
and excluding the eigenvalue resulting from the removed vertex ($0$ for a
duplicate or $-1$ for a co-duplicate).
The displacement is also constrained by the interlacing theorem.
The two roots of \eqref{eqn:quadratic} are either both real or else they
are complex conjugates.
In the case of real roots, we first exclude roots that lie outside the range
allowed by the interlacing theorem.
If both roots are within the allowed range, the value closest to the first
order approximation is taken as the estimate.
For complex conjugate roots, the real part is taken instead.
The easily obtained values $f(\lambda_0)$, $f'(\lambda_0)$, and
$f''(\lambda_0)$ allow us to obtain an estimate for the eigenvalues
of $G_{-v_\ell}$ without solving the high-order polynomial equation
$f(\lambda)=0$.



\appendix

\section{Examples on Nested Split Graphs}
\label{sec:examples_nsg}

We illustrate the use of the results from
Section~\ref{sec:displacement} on examples from the class of nested
split graphs (NSG), also known in the literature as threshold graphs.
Following the notation of \cite{sbd18dam}, the compact creation sequence
is $\mathbf{a} = (a_1,a_2,\ldots,a_r)$, where $\sum a_i=n$, the number of
cells $r$ is even, and $a_i \geq 1 \ \forall i$.
This represents the connected graph
$( \cdots
 (( \overline{K_{a_1}} \bigtriangledown K_{a_2} )
 \dotcup \overline{K_{a_3}} )
 \cdots
 \dotcup \overline{K_{a_{r-1}}} ) \bigtriangledown K_{a_r}$
where $K_s$ is the complete graph on $s$ vertices, $\overline{K_s}$ is
its complement, while $\bigtriangledown$ and $\dotcup$ are the graph
operators join and disjoint union respectively.
Note that $\mathbf{a}$ has $r$ cells, of which $(a_1,a_3,\ldots,a_{r-1})$
are co-clique cells and $(a_2,a_4,\ldots,a_r)$ are clique cells.
A NSG with $r$ cells has $r$ main eigenvalues if $a_1 \geq 2$
and $r-1$ if $a_1 = 1$.
Recall that a \emph{main} eigenvalue of a graph $G$ is an eigenvalue $\mu$
of $\vec{A}$ such that $\vec{A}$ has some eigenvector $\vec{x}$ not orthogonal
to the all-one vector $\vec{j}$ associated with $\mu$ \cite{rowlinson2007main,
sciriha2012spectrum}.
The significance of the non-zero main eigenvalues is that they determine
the number of walks of any length in $G$ \cite{cvetkovic2001introduction,
harary1979spectral}.
In a NSG, the spectrum consists of the main eigenvalues (except $0$ or
$-1$, which are never main in a NSG), the eigenvalue zero with multiplicity
determined by the duplicate vertices, and the eigenvalue $-1$ with multiplicity
determined by the co-duplicate vertices.

The following examples consider different operations on the NSG $G$, having 18
vertices in 10 cells, with compact creation sequence
$\vec{a} = (2, 2, 2, 2, 2, 2, 2, 2, 1, 1)$.
This graph therefore has 10 main eigenvalues.
Its characteristic polynomial is
\begin{align*}
\cp{\vec{A}(G)}
&= \lambda^{4} (\lambda+1)^4 \left(
\lambda^{10} - 4 \lambda^{9} - 75 \lambda^{8} - 128 \lambda^{7} + 371 \lambda^{6} + 860 \lambda^{5} - 441 \lambda^{4} \right.\\
&\quad\left.
- 1368 \lambda^{3} + 336 \lambda^{2} + 704 \lambda - 256
\right)
\text{.}
\end{align*}

\subsection{Removing a duplicate vertex}
\label{sec:example1}

Consider deleting a vertex from the third cell of the graph $G$, resulting in a
graph $G'$ with compact creation sequence given by
$\vec{a'} = (2, 2, 1, 2, 2, 2, 2, 2, 1, 1)$, having 17 vertices in 10 cells.
When listing the vertices in the same order in the adjacency matrix, this
means that we are removing one of vertices 5 or 6, which are duplicates.
From Theorem~\ref{thm:Glk} we can obtain the characteristic polynomial of $G'$
from that of $G$ by first dividing by $\lambda$ to remove a zero eigenvalue,
then adding $h_{5,6}$ to obtain the necessary displacement in the remaining
eigenvalues.
Using Proposition~\ref{prop:hlk} we obtain
\begin{align*}
h_{5,6}
&= 7 \lambda^{15} + 42 \lambda^{14} + 20 \lambda^{13} - 348 \lambda^{12} - 758 \lambda^{11} + 192 \lambda^{10} + 2220 \lambda^{9} + 2124 \lambda^{8} \nonumber\\
&\quad
- 489 \lambda^{7} - 1722 \lambda^{6} - 616 \lambda^{5} + 224 \lambda^{4} + 128 \lambda^{3}
\text{.}
\end{align*}
It can be verified that applying Theorem~\ref{thm:Glk} gives
\begin{align*}
\cp{\vec{A}(G')}
&= \lambda^{3} (\lambda+1)^4 \left(
\lambda^{10} - 4 \lambda^{9} - 68 \lambda^{8} - 114 \lambda^{7} + 293 \lambda^{6} + 712 \lambda^{5} - 202 \lambda^{4} \right.\\
&\quad\left.
- 946 \lambda^{3} + 104 \lambda^{2} + 416 \lambda - 128
\right)
\text{.}
\end{align*}
We can now estimate the shift in the main eigenvalues from $G$ to $G'$ using
the method of Section~\ref{sec:displacement}, after obtaining the necessary
functions $f(\lambda)$, $f'(\lambda)$, and $f''(\lambda)$.
Table~\ref{tab:example1} gives the main eigenvalues of $G$ and $G'$, the actual
displacement, and the estimates computed using the first-order and second-order
approximations of Section~\ref{sec:displacement}.
\inserttab{tab:example1}{ccccc}
{Removing a duplicate vertex:
the eigenvalues of $G$ with compact creation sequence
$\vec{a} = (2, 2, 2, 2, 2, 2, 2, 2, 1, 1)$
and $G'$ with
$\vec{a'} = (2, 2, 1, 2, 2, 2, 2, 2, 1, 1)$,
the actual displacement, and estimates using first- and second-order
approximations.}{
\multicolumn{2}{c}{Eigenvalues} & Actual & \multicolumn{2}{c}{Estimates} \\
$G$ & $G'$ & Displacement & First-order & Second-order%
\footnote{The chosen estimate is shown in bold.} \\
\hline
$-4.45$ & $-4.05$ & $0.398$     & $0.151$     & $\mathbf{0.182} \pm 0.148 j$    \\
$-2.28$ & $-2.09$ & $0.182$     & $0.0671$    & $\mathbf{0.0819} \pm 0.0655 j$  \\
$-1.76$ & $-1.72$ & $0.0377$    & $0.0265$    & $0.0660, \mathbf{0.0443}$       \\
$-1.5$  & $-1.43$ & $0.0673$    & $0.0304$    & $\mathbf{0.0502} \pm 0.0231 j$  \\
$-1.43$ & $-1.42$ & $0.00937$   & $-0.00148$  & $\mathbf{0.00880}, -0.00127$    \\

$-1$ & $-1$%
\footnote{Repeated 4 times in $G$ and $G'$, comparison unnecessary.}
& $0$ & --- & --- \\

$0$ & $0$%
\footnote{Repeated 4 times in $G$, 3 times in $G'$, comparison unnecessary.}
& $0$ & --- & --- \\

$0.432$ & $0.431$ & $-0.000233$ & $-0.000233$ & $0.419, \mathbf{-0.000233}$     \\
$0.697$ & $0.52$  & $-0.178$    & $-0.0823$   & $-0.223, \mathbf{-0.131}$       \\
$1$     & $0.901$ & $-0.0990$   & $-0.0513$   & $\mathbf{-0.0736} \pm 0.0462 j$ \\
$1.96$  & $1.95$  & $-0.0116$   & $-0.0108$   & $-0.152, \mathbf{-0.0117}$      \\
$11.3$  & $10.9$  & $-0.406$    & $-0.262$    & $\mathbf{-0.456} \pm 0.176 j$   \\
}

\subsection{A special case of removing a duplicate vertex}
\label{sec:example2}

In this example, we delete a vertex from the first cell of the graph $G$,
resulting in a graph $G'$ with compact creation sequence given by
$\vec{a'} = (1, 2, 2, 2, 2, 2, 2, 2, 1, 1)$.
This is a special case, because a single vertex in the first cell is
effectively a co-duplicate of the vertices in the second cell.
As a result, the number of main eigenvalues decreases by one.
As in the general case, we obtain the characteristic polynomial of $G'$
from that of $G$ by first dividing by $\lambda$ to remove a zero eigenvalue,
then adding $h_{1,2}$ to obtain the necessary displacement in the remaining
eigenvalues.
In this case, however, we know that the number of main eigenvalues reduces
by one and the number of eigenvalues $-1$ increases by one, as effectively
an additional co-duplicate is created.
That is, we do not need one of the estimates that will be calculated.
So, proceeding as in the earlier example, using Proposition~\ref{prop:hlk}
we obtain
\begin{align*}
h_{1,2}
&= 9 \lambda^{15} + 72 \lambda^{14} + 140 \lambda^{13} - 280 \lambda^{12} - 1370 \lambda^{11} - 1304 \lambda^{10} + 1228 \lambda^{9} + 2840 \lambda^{8} \nonumber\\
&\quad
+ 793 \lambda^{7} - 1328 \lambda^{6} - 800 \lambda^{5} + 128 \lambda^{4} + 128 \lambda^{3}
\text{.}
\end{align*}
Using Theorem~\ref{thm:Glk} this gives
\begin{align*}
\cp{\vec{A}(G')}
&= \lambda^{3} (\lambda+1)^5 \left(
\lambda^{9} - 5 \lambda^{8} - 61 \lambda^{7} - 31 \lambda^{6} + 344 \lambda^{5} + 216 \lambda^{4} - 632 \lambda^{3} \right.\nonumber\\
&\quad\left.
- 144 \lambda^{2} + 448 \lambda - 128
\right)
\text{.}
\end{align*}
Estimates for the shift in the main eigenvalues from $G$ to $G'$ using the
first-order and second-order approximations of Section~\ref{sec:displacement}
are given in Table~\ref{tab:example2}, together with the main eigenvalues
of $G$ and $G'$ and the actual displacement.
\inserttab{tab:example2}{ccccc}
{Removing a duplicate vertex -- special case:
the eigenvalues of $G$ with compact creation sequence
$\vec{a} = (2, 2, 2, 2, 2, 2, 2, 2, 1, 1)$
and $G'$ with
$\vec{a'} = (1, 2, 2, 2, 2, 2, 2, 2, 1, 1)$,
the actual displacement, and estimates using first- and second-order
approximations.}{
\multicolumn{2}{c}{Eigenvalues} & Actual & \multicolumn{2}{c}{Estimates} \\
$G$ & $G'$ & Displacement & First-order & Second-order%
\footnote{The chosen estimate is shown in bold.} \\
\hline
$-4.45$ & $-4.24$ & $0.209$                & $0.113$                & $\mathbf{0.162} \pm 0.101 j$            \\
$-2.28$ & $-2.2$  & $0.0766$               & $0.0452$               & $\mathbf{0.0713} \pm 0.0369 j$          \\
$-1.76$ & $-1.73$ & $0.0275$               & $0.0214$               & $0.0840, \mathbf{0.0288}$               \\
$-1.5$  & $-1.43$ & $0.0686$               & $0.0450$               & $0.110, \mathbf{0.0759}$                \\
$-1.43$ & $-1$%
\footnote{Multiplicity of $-1$ is known to increase by one, comparison unnecessary.}
& $0.432$ & --- & --- \\

$-1$ & $-1$%
\footnote{Repeated 4 times in $G$ and $G'$, comparison unnecessary.}
& $0$ & --- & --- \\

$0$ & $0$%
\footnote{Repeated 4 times in $G$, 3 times in $G'$, comparison unnecessary.}
& $0$ & --- & --- \\

$0.432$ & $0.432$ & $-2.03 \times 10^{-6}$ & $-2.03 \times 10^{-6}$ & $-0.188, \mathbf{-2.03 \times 10^{-6}}$ \\
$0.697$ & $0.683$ & $-0.0145$              & $-0.0131$              & $-0.130, \mathbf{-0.0145}$              \\
$1$     & $0.951$ & $-0.0486$              & $-0.0323$              & $\mathbf{-0.0598} \pm 0.0167 j$         \\
$1.96$  & $1.85$  & $-0.116$               & $-0.0663$              & $\mathbf{-0.101} \pm 0.0568 j$          \\
$11.3$  & $10.7$  & $-0.634$               & $-0.341$               & $\mathbf{-0.490} \pm 0.307 j$           \\
}

\subsection{Removing a co-duplicate vertex}
\label{sec:example3}

Finally, we delete a vertex from the second cell of the graph $G$, resulting in
a graph $G'$ with compact creation sequence given by
$\vec{a'} = (2, 1, 2, 2, 2, 2, 2, 2, 1, 1)$.
In this case we are removing a co-duplicate vertex, so we obtain the
characteristic polynomial of $G'$ from that of $G$ by first dividing by
$\lambda+1$ to remove one of the $-1$ eigenvalues, then adding $h_{3,4}$
to obtain the necessary displacement in the remaining eigenvalues.
So, proceeding as in the earlier examples, using Proposition~\ref{prop:hlk}
we obtain
\begin{align*}
h_{3,4}
&= \lambda^{16} + 9 \lambda^{15} + 4 \lambda^{14} - 171 \lambda^{13} - 596 \lambda^{12} - 507 \lambda^{11} + 888 \lambda^{10} + 1923 \lambda^{9} + 599 \lambda^{8} \nonumber\\
&\quad
- 1062 \lambda^{7} - 736 \lambda^{6} + 96 \lambda^{5} + 128 \lambda^{4}
\text{.}
\end{align*}
Using Theorem~\ref{thm:Glk} this gives
\begin{align*}
\cp{\vec{A}(G')}
&= \lambda^{4} (\lambda+1)^3 \left(
\lambda^{10} - 3 \lambda^{9} - 69 \lambda^{8} - 145 \lambda^{7} + 232 \lambda^{6} + 726 \lambda^{5} - 112 \lambda^{4} \right.\nonumber\\
&\quad\left.
- 926 \lambda^{3} + 80 \lambda^{2} + 416 \lambda - 128
\right)
\text{.}
\end{align*}
Estimates for the shift in the main eigenvalues from $G$ to $G'$ using the
first-order and second-order approximations of Section~\ref{sec:displacement}
are given in Table~\ref{tab:example3}, together with the main eigenvalues
of $G$ and $G'$ and the actual displacement.
\inserttab{tab:example3}{ccccc}
{Removing a co-duplicate vertex:
the eigenvalues of $G$ with compact creation sequence
$\vec{a} = (2, 2, 2, 2, 2, 2, 2, 2, 1, 1)$
and $G'$ with
$\vec{a'} = (2, 1, 2, 2, 2, 2, 2, 2, 1, 1)$,
the actual displacement, and estimates using first- and second-order
approximations.}{
\multicolumn{2}{c}{Eigenvalues} & Actual & \multicolumn{2}{c}{Estimates} \\
$G$ & $G'$ & Displacement & First-order & Second-order%
\footnote{The chosen estimate is shown in bold.} \\
\hline
$-4.45$ & $-4.34$ & $0.101$                & $0.0716$               & $0.162, \mathbf{0.128}$                 \\
$-2.28$ & $-2.27$ & $0.00308$              & $0.00299$              & $0.108, \mathbf{0.00308}$               \\
$-1.76$ & $-1.76$ & $0.00206$              & $0.00201$              & $0.0899, \mathbf{0.00206}$              \\
$-1.5$  & $-1.43$ & $0.0685$               & $0.0382$               & $\mathbf{0.0659} \pm 0.0262 j$          \\
$-1.43$ & $-1.35$ & $0.0818$               & $-7.45 \times 10^{-5}$ & $\mathbf{0.0520}, -7.44 \times 10^{-5}$ \\

$-1$ & $-1$%
\footnote{Repeated 4 times in $G$, 3 times in $G'$, comparison unnecessary.}
& $0$ & --- & --- \\

$0$ & $0$%
\footnote{Repeated 4 times in $G$ and $G'$, comparison unnecessary.}
& $0$ & --- & --- \\

$0.432$ & $0.432$ & $-3.73 \times 10^{-5}$ & $-3.73 \times 10^{-5}$ & $-0.353, \mathbf{-3.73 \times 10^{-5}}$ \\
$0.697$ & $0.567$ & $-0.131$               & $-0.0759$              & $-0.474, \mathbf{-0.0903}$              \\
$1$     & $0.85$  & $-0.150$               & $-0.0608$              & $\mathbf{-0.0781} \pm 0.0583 j$         \\
$1.96$  & $1.74$  & $-0.227$               & $-0.0921$              & $\mathbf{-0.114} \pm 0.0894 j$          \\
$11.3$  & $10.6$  & $-0.748$               & $-0.372$               & $\mathbf{-0.504} \pm 0.347 j$           \\
}

It may come as a surprise that there are shifts in most of the eigenvalues
when removing a twin vertex.
Considering the limited interval in which the maximum eigenvalue can lie,
we note that its displacement when the graph is perturbed is significant.


\section{Examples on General Graphs}
\label{sec:examples_nut}

\begin{figure*}[!tp]
   \centering
   \includegraphics[width=0.5\textwidth]{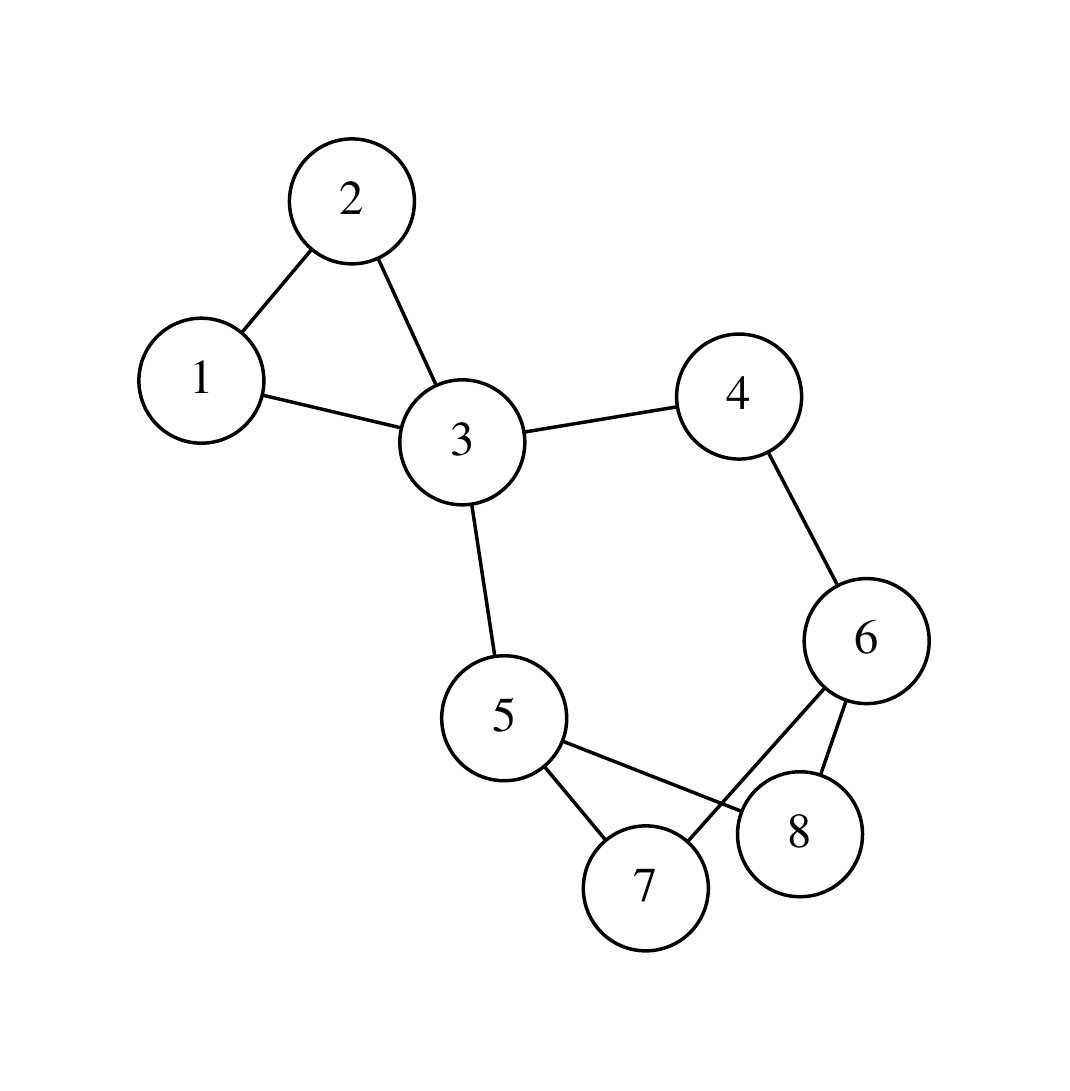}
   \caption{A general graph $G$ with two duplicates (vertices 7 and 8) and two
   co-duplicates (vertices 1 and 2).}
   \label{fig:G8}
\end{figure*}

We also illustrate the use of the results from Section~\ref{sec:displacement}
on a more general graph $G$, shown in Figure~\ref{fig:G8}.
This graph has 6 main eigenvalues, and its characteristic polynomial is
\begin{align*}
\cp{\vec{A}(G)}
&= \lambda (\lambda+1) \left(
\lambda^{6} - \lambda^{5} - 9 \lambda^{4} + 7 \lambda^{3} + 19 \lambda^{2} - 13 \lambda
\right)
\text{.}
\end{align*}

\subsection{Removing a co-duplicate vertex}
\label{sec:example4}

Consider deleting a co-duplicate vertex (1 or 2) from graph $G$, resulting
in a graph $G'$.
We obtain the characteristic polynomial of $G'$ from that of $G$ by first
dividing by $\lambda+1$ to remove a $-1$ eigenvalue, then adding $h_{1,2}$
to obtain the necessary displacement in the remaining eigenvalues.
So, proceeding as in the earlier examples, using Proposition~\ref{prop:hlk}
we obtain
\begin{align*}
h_{1,2}
&= \lambda^{6} + \lambda^{5} - 7 \lambda^{4} - 5 \lambda^{3} + 9 \lambda^{2} - 2 \lambda
\text{.}
\end{align*}
Using Theorem~\ref{thm:Glk} this gives
\begin{align*}
\cp{\vec{A}(G')}
&= \lambda \left(
\lambda^{6} - 8 \lambda^{4} + 14 \lambda^{2} - 4 \lambda - 2
\right)
\text{.}
\end{align*}
We can now estimate the shift in the main eigenvalues from $G$ to $G'$ using
the method of Section~\ref{sec:displacement}, after obtaining the necessary
functions $f(\lambda)$, $f'(\lambda)$, and $f''(\lambda)$.
Table~\ref{tab:example4} gives the main eigenvalues of $G$ and $G'$, the actual
displacement, and the estimates computed using the first-order and second-order
approximations of Section~\ref{sec:displacement}.
\inserttab{tab:example4}{ccccc}
{Removing a co-duplicate vertex:
the eigenvalues of graph $G$ of Figure~\ref{fig:G8}
and $G'$ obtained by removing vertex 1 or 2,
the actual displacement, and estimates using first- and second-order
approximations.}{
\multicolumn{2}{c}{Eigenvalues} & Actual & \multicolumn{2}{c}{Estimates} \\
$G$ & $G'$ & Displacement & First-order & Second-order%
\footnote{The chosen estimate is shown in bold.} \\
\hline
$-2.20$ & $-2.18$ & $0.0139$   & $0.0130$   & $0.214, \mathbf{0.0139}$    \\
$-1.89$ & $-1.83$ & $0.0592$   & $0.0642$   & $-0.613, \mathbf{0.0581}$  \\

$-1$%
\footnote{Due to co-duplicate in $G$; removed in $G'$.}
& --- & --- & --- & --- \\

$0$     & $-0.265$ & $-0.265$  & $0$        & -0.5, $\mathbf{0}$  \\

$0$ & $0$%
\footnote{Due to duplicate in $G$; remains in $G'$.}
& $0$ & --- & --- \\

$0.664$ & $0.656$ & $-0.00768$ & $-0.00763$ & $-1.26, \mathbf{-0.00768}$     \\
$1.79$  & $1.18$  & $-0.609$   & $-0.385$   & $-17.6, \mathbf{-0.394}$       \\
$2.64$  & $2.45$  & $-0.192$   & $-0.132$   & $\mathbf{-0.262} \pm 0.0261 j$      \\
}

\subsection{Removing a duplicate vertex}
\label{sec:example5}

Finally, we delete a duplicate vertex (7 or 8) from graph $G$, resulting in
a graph $G'$.
In this case we obtain the characteristic polynomial of $G'$ from that of
$G$ by first dividing by $\lambda$ to remove one of the $0$ eigenvalues,
then adding $h_{7,8}$ to obtain the necessary displacement in the remaining
eigenvalues.
So, proceeding as in the earlier examples, using Proposition~\ref{prop:hlk}
we obtain
\begin{align*}
h_{7,8}
&= 2 \lambda^{5} - 7 \lambda^{3} - 4 \lambda^{2} + 3 \lambda + 2
\text{.}
\end{align*}
Using Theorem~\ref{thm:Glk} this gives
\begin{align*}
\cp{\vec{A}(G')}
&= (\lambda+1) \left(
\lambda^{6} - \lambda^{5} - 7 \lambda^{4} + 5 \lambda^{3} + 11 \lambda^{2} - 7 \lambda
\right)
\text{.}
\end{align*}
Estimates for the shift in the main eigenvalues from $G$ to $G'$ using the
first-order and second-order approximations of Section~\ref{sec:displacement}
are given in Table~\ref{tab:example5}, together with the main eigenvalues
of $G$ and $G'$ and the actual displacement.
\inserttab{tab:example5}{ccccc}
{Removing a duplicate vertex:
the eigenvalues of graph $G$ of Figure~\ref{fig:G8}
and $G'$ obtained by removing vertex 7 or 8,
the actual displacement, and estimates using first- and second-order
approximations.}{
\multicolumn{2}{c}{Eigenvalues} & Actual & \multicolumn{2}{c}{Estimates} \\
$G$ & $G'$ & Displacement & First-order & Second-order%
\footnote{The chosen estimate is shown in bold.} \\
\hline
$-2.20$ & $-1.94$ & $0.262$   & $0.265$   & $\mathbf{0.241} \pm 0.264 j$    \\
$-1.89$ & $-1.62$ & $0.388$   & $0.388$   & $\mathbf{0.212} \pm 0.346 j$  \\

$-1$ & $-1$%
\footnote{Due to co-duplicate in $G$; remains in $G'$.}
& --- & --- & --- \\

$0$%
\footnote{Due to duplicate in $G$; removed in $G'$.}
& --- & --- & --- & --- \\

$0$     & $0$     & $0$%
\footnote{Displacement constrained by interlacing; no estimate required.}
& --- & --- \\
$0.664$ & $0.618$ & $-0.0458$ & $-0.0531$ & $-0.538, \mathbf{-0.0589}$     \\
$1.79$  & $1.46$  & $-0.322$  & $-0.193$  & $-5.01, \mathbf{-0.201}$       \\
$2.64$  & $2.47$  & $-0.167$  & $-0.162$  & $\mathbf{-0.263} \pm 0.127 j$      \\
}



\bibliographystyle{elsarticle-num}
\bibliography{references}

\begin{thebibliography}{10}
\expandafter\ifx\csname url\endcsname\relax
  \def\url#1{\texttt{#1}}\fi
\expandafter\ifx\csname urlprefix\endcsname\relax\def\urlprefix{URL }\fi
\expandafter\ifx\csname href\endcsname\relax
  \def\href#1#2{#2} \def\path#1{#1}\fi

\bibitem{heilbronner1953kompositions}
E.~Heilbronner, Das kompositions-prinzip: Eine anschauliche methode zur
  elektronen-theoretischen behandlung nicht oder niedrig symmetrischer molekeln
  im rahmen der mo-theorie, Helvetica Chimica Acta 36~(1) (1953) 170--188.

\bibitem{heilbronner1954molecular}
E.~Heilbronner, Molecular orbitals in homologen reihen mehrkerniger
  aromatischer kohlenwasserstoffe: I. die eigenwerte von {LCAO-MO}'s in
  homologen reihen, Helvetica Chimica Acta 37~(3) (1954) 921--935.

\bibitem{heilbronner1954graphisches}
E.~Heilbronner, Ein graphisches verfahren zur faktorisierung der
  s{\"a}kulardeterminante aromatischer ringsysteme im rahmen der
  {LCAO--MO}-theorie, Helvetica Chimica Acta 37~(3) (1954) 913--921.

\bibitem{heilbronner1962graphentheoretischen}
E.~Heilbronner, {\"U}ber einen graphentheoretischen zusammenhang zwischen dem
  h{\"u}ckel'schen {MO}-verfahren und dem formalismus der resonanztheorie,
  Helvetica Chimica Acta 45~(5) (1962) 1722--1725.

\bibitem{heilbronner1979some}
E.~Heilbronner, Some comments on cospectral graphs, Math. Chem 5 (1979)
  105--133.

\bibitem{schwenk1974}
A.~J. Schwenk, Computing the characteristic polynomial of a graph, in: Graphs
  and combinatorics, Springer, 1974, pp. 153--172.

\bibitem{rosenfeld2013}
V.~R. Rosenfeld, Another form of the transmission function, Journal of
  Mathematical Chemistry 51~(10) (2013) 2639--2643.

\bibitem{thune1979}
M.~Th{\"u}ne, Eigenvalues of matrices and graphs, Ph.D. thesis (1979).

\bibitem{haemers1995}
W.~H. Haemers, Interlacing eigenvalues and graphs, Linear Algebra and its
  applications 226 (1995) 593--616.

\bibitem{haemers1980}
W.~H. Haemers, Eigenvalue techniques in design and graph theory, Ph.D. thesis
  (1980).

\bibitem{lovasz2007}
L.~Lov{\'a}sz, \href{http://web.cs.elte.hu/~lovasz/eigenvals-x.pdf}{Eigenvalues
  of graphs} (2007).
\newline\urlprefix\url{http://web.cs.elte.hu/~lovasz/eigenvals-x.pdf}

\bibitem{marino2006}
M.~C. Marino, I.~Sciriha, S.~K. Simi{\'c}, D.~V. To\v{s}i{\'c}, More about
  singular line graphs of trees, Publications de l'Institut Mathématique
  79~(93).

\bibitem{sdbfp2013}
I.~Sciriha, M.~Debono, M.~Borg, P.~Fowler, B.~Pickup, Interlacing–-extremal
  graphs, Ars Math. Contemp. 6 (2013) 261--278.
\newblock \href {http://dx.doi.org/10.26493/1855-3974.275.574}
  {\path{doi:10.26493/1855-3974.275.574}}.

\bibitem{rowlinson2007}
P.~Rowlinson, Co-cliques and star complements in extremal strongly regular
  graphs, Linear algebra and its applications 421~(1) (2007) 157--162.

\bibitem{sciriha2012spectrum}
I.~Sciriha, S.~Farrugia, On the spectrum of threshold graphs, ISRN Discrete
  Mathematics 2011.

\bibitem{mohammadian2016}
A.~Mohammadian, V.~Trevisan, Some spectral properties of cographs, Discrete
  Mathematics 339~(4) (2016) 1261--1264.

\bibitem{so1999}
W.~So, Rank one perturbation and its application to the laplacian spectrum of a
  graph, Linear and Multilinear Algebra 46~(3) (1999) 193--198.

\bibitem{froberg}
C.-E. Fr{\"o}berg, Introduction to numerical analysis, Vol.~6, Addison-Wesley
  Reading, MA, 1965.

\bibitem{sbd18dam}
I.~Sciriha, J.~A. Briffa, M.~Debono, Fast algorithms for indices of nested
  split graphs approximating real complex networks, Discrete Applied
  Mathematics 247 (2018) 152--164.
\newblock \href {http://dx.doi.org/10.1016/j.dam.2018.03.054}
  {\path{doi:10.1016/j.dam.2018.03.054}}.

\bibitem{rowlinson2007main}
P.~Rowlinson, The main eigenvalues of a graph: a survey, Applicable Analysis
  and Discrete Mathematics (2007) 455--471.

\bibitem{cvetkovic2001introduction}
D.~Cvetkovi{\'c}, P.~Rowlinson, S.~Simi{\'c}, An Introduction to the Theory of
  Graph Spectra (London Mathematical Society Student Texts), Cambridge:
  Cambridge University Press, 2001.

\bibitem{harary1979spectral}
F.~Harary, A.~Schwenk, The spectral approach to determining the number of walks
  in a graph, Pacific Journal of Mathematics 80~(2) (1979) 443--449.

\end{thebibliography}

\end{document}